\documentclass[12pt]{amsart}

\usepackage{amssymb}

\newtheorem{theorem}{Theorem}[section]
\newtheorem{proposition}[theorem]{Proposition}
\newtheorem{lemma}[theorem]{Lemma}
\newtheorem{corollary}[theorem]{Corollary}
\newtheorem{remark}[theorem]{Remark}

\numberwithin{equation}{section}

\begin{document}
\baselineskip=15pt

\title[Semistable bundle over projective manifolds]{On
semistable principal bundles\\[5pt] over a complex projective manifold}

\author[I. Biswas]{Indranil Biswas}

\address{School of Mathematics, Tata Institute of Fundamental
Research, Homi Bhabha Road, Bombay 400005, India}

\email{indranil@math.tifr.res.in}

\author[U. Bruzzo]{Ugo Bruzzo}

\address{Scuola Internazionale Superiore di Studi Avanzati,
Via Beirut 2--4, 34013, Trieste, Italy}

\email{bruzzo@sissa.it}

\subjclass[2000]{32L05, 14L10, 14F05}

\keywords{Principal bundle, semistability, numerically effectiveness}

\date{}

\begin{abstract}

Let $G$ be a simple linear algebraic group defined over
the field of complex numbers. Fix a proper parabolic subgroup $P$
of $G$, and also fix a nontrivial
antidominant character $\chi$ of $P$. We prove that a holomorphic
principal $G$--bundle $E_G$ over a
connected complex projective manifold $M$ is semistable
satisfying the condition that the second Chern class
$c_2(\text{ad}(E_G))\, \in \, H^4(M,\, {\mathbb Q})$ vanishes
if and only if the line bundle over $E_G/P$ defined by
$\chi$ is numerically effective. Also, a
principal $G$--bundle $E_G$ over $M$
is semistable with $c_2(\text{ad}(E_G))\, =\, 0$ if and only if for
every pair of the form
$(Y\, ,\psi)$, where $\psi$ is a holomorphic map
to $M$ from a compact connected Riemann surface $Y$,
and for every holomorphic reduction of structure group
$E_P\, \subset\, \psi^*E_G$ to the subgroup $P$,
the line bundle over $Y$ associated to
the principal $P$--bundle $E_P$ for $\chi$
is of nonnegative degree. Therefore, $E_G$
is semistable with $c_2(\text{ad}(E_G))\, =\, 0$ if and
only if for each pair $(Y\, ,\psi)$ of the above type the
$G$--bundle $\psi^*E_G$ over $Y$ is semistable.

Similar results remain valid for principal
bundles over $M$ with a reductive linear
algebraic group as the
structure group. These generalize an earlier work of
Y. Miyaoka, \cite{Mi}, where he gave a characterization
of semistable vector bundles over a smooth projective curve.
Using these characterizations one can also
produce similar criteria for the semistability of parabolic
principal bundles over a compact Riemann surface.
\end{abstract}

\maketitle

\section{Introduction}

Let $G$ be a simple linear algebraic group
defined over $\mathbb C$.
Fix a proper parabolic subgroup $P$ of $G$. Fix a nontrivial
antidominant character $\chi$ of $P$. This means that
the associated line bundle $(G\times {\mathbb C}_\chi)/P$ over
$G/P$ is numerically effective and nontrivial.

Let $M$ be a connected smooth complex projective variety. Let
$E_G$ be a holomorphic principal $G$--bundle over $M$.
The natural projection $E_G\, \longrightarrow\, E_G/P$
defines a principal $P$--bundle over $E_G/P$, and hence the
character $\chi$ of $P$ associates a holomorphic line bundle
$$
L_\chi\,=\, (E_G\times {\mathbb C})/P
$$
over $E_G/P$, where the action of $p\,\in\, P$ sends any $(z\, ,
c) \,\in\, E_G\times {\mathbb C}$ to $(zp\, ,\chi(p)^{-1}c)$.

We prove the following criterion for $L_\chi$ to be
numerically effective (see Theorem \ref{thm.h.d.}):

\begin{theorem}
The line bundle $L_\chi$ over $E_G/P$
is numerically effective if and only if $E_G$ is semistable and
$c_2(\text{ad}(E_G))\, =\, 0$.
\end{theorem}

It should be clarified that this criterion for semistability
of any $G$--bundle $E_G$ with $c_2(\text{ad}(E_G))\, =\, 0$
needs to be verified for just one pair $(P\, ,\chi)$.
Since the condition $c_2(\text{ad}(E_G))\, =\, 0$ is automatically
satisfied when $M$ is a Riemann surface,
Theorem \ref{thm.h.d.} gives a characterization of semistable
$G$--bundles over a compact connected Riemann surface.

In the special case where
$\dim M \, =\, 1$, $G\,=\, \text{GL}(n,{\mathbb C})$, and $P$
is the parabolic subgroup of $\text{GL}(n,{\mathbb C})$ that
fixes a line in ${\mathbb C}^n$ (so $G/P\,=\,
{\mathbb C}{\mathbb P}^{n-1}$), Theorem \ref{thm.h.d.}
is due to Y. Miyaoka (see \cite[page 456, Theorem 3.1]{Mi});
in this case $\text{Pic}(G/P)\, =\, {\mathbb Z}$ and hence
$\chi$ is unique up to taking tensor powers. In \cite{BH}, Miyaoka's
result was generalized to the context of Higgs vector bundles
over curves.

As a consequence of Theorem \ref{thm.h.d.}, the line bundle
$L_\chi$ over $E_G/P$ is
numerically effective if and only if for every parabolic
subgroup $P'\, \subset\, G$, and every nontrivial antidominant
character $\chi'$ of $P'$, the associated line bundle
$L_{\chi'}\, :=\, (E_G\times{\mathbb C}_{\chi'})/P'$ over
$E_G/P'$ is numerically effective. The above condition that
$c_2(\text{ad}(E_G))\, \in\, H^4(M,\, {\mathbb Q})$
vanishes is equivalent to the condition that the real
characteristic class of $E_G$ corresponding to the Killing
form on the Lie algebra of $G$ vanishes (Remark \ref{rem.-2}).

In Section \ref{s.12} we prove Theorem \ref{thm.h.d.}
under the assumption that $\dim M \, =\, 1$ (see Proposition
\ref{pro.11}). The proof of Theorem \ref{thm.h.d.}
given in Section \ref{s.13} crucially uses this special case.

Theorem \ref{thm.h.d.} can be reformulated as follows:

A $G$--bundle $E_G$ over $M$ is semistable with
$c_2(\text{ad}(E_G))\, =\, 0$ if and only if
for every pair of the form $(Y\, ,\psi)$, where $Y$ is a
compact connected Riemann surface and
$$
\psi\,:\, Y\,\longrightarrow\, M
$$
a holomorphic map, and for every reduction
$E_P\, \subset\, \psi^*E_G$ of structure group, to
the fixed parabolic subgroup $P$, of the principal
$G$--bundle $\psi^*E_G$ over $Y$, the associated line bundle
$E_P(\chi)\, =\, (E_P\times {\mathbb C})/P$ over $Y$ is of
nonnegative degree, where $\chi$ is the
fixed character of $P$ (Proposition \ref{pr.-1}).

It is known that a $G$--bundle $E_G$ over $M$ is semistable
if and only if the restriction of $E_G$ to the general complete
intersection curve of sufficiently large degree hypersurfaces
is semistable. Therefore, we may ask the following question:
under what condition the restriction of $E_G$ to every
smooth curve in $M$ is semistable? As a
consequence of Proposition \ref{pr.-1} and Theorem
\ref{thm.h.d.} this question has the following answer:

\textit{A $G$--bundle $E_G$ is semistable and
$c_2({\rm ad}(E_G))\, =\, 0$ if and only if for each
pair $(Y\, ,\psi)$ of the above type the $G$--bundle
$\psi^*E_G$ over $Y$ is semistable.}

In Section \ref{sec.-4} we give
the following analog of Theorem
\ref{thm.h.d.} and Proposition \ref{pr.-1} for principal
bundles with a reductive group as the structure group
(see Theorem \ref{thm.2}):

\begin{theorem}
Let $G$ be a connected reductive linear algebraic group over
$\mathbb C$. Fix a parabolic subgroup $P\, \subset\, G$ without
any simple factor, and also fix a character $\chi$
of $P$ such that
\begin{enumerate}
\item[(i)] $\chi$ is trivial on the center $Z(G)\, \subset\, G$,
and
\item[(ii)] the restriction of $\chi$ to the parabolic subgroup
of each simple factor of $G/Z(G)$ defined by $P$ is
nontrivial and antidominant.
\end{enumerate}
Let $E_G$ be a principal $G$--bundle over a connected projective
manifold $M$.
Then the following four statements are equivalent:
\begin{enumerate}
\item{} The $G$--bundle $E_G$ is semistable and the
second Chern class
$$
c_2({\rm ad}(E_G)) \, \in\, H^4(M,\, {\mathbb Q})
$$
vanishes.

\item{} The associated line bundle $L_\chi\, :=\,
(E_G\times {\mathbb C}_{\chi})/P$ over $E_G/P$ for the
character $\chi$ is numerically effective.

\item{} For every pair of the form $(Y\, ,\psi)$,
where $Y$ is a compact connected Riemann surface and
$$
\psi\,:\, Y\,\longrightarrow\, M
$$
a holomorphic map, and every reduction $E_P\, \subset\,
\psi^*E_G$ of structure group to $P$ of the principal
$G$--bundle $\psi^*E_G$ over $Y$, the associated line bundle
$E_P(\chi)\, =\, (E_P\times {\mathbb C}_{\chi})/P$ over
$Y$ is of nonnegative degree.

\item{} For any pair $(Y\, ,\psi)$ as in (3), the $G$--bundle
$\psi^*E_G$ over $Y$ is semistable.
\end{enumerate}
\end{theorem}

Let $G$ be a connected reductive linear algebraic group
over $\mathbb C$. A \emph{parabolic $G$-bundle} $E_\ast$
over a connected smooth complex projective curve $X$
with parabolic structure over a reduced
divisor $D\subset X$ is, loosely speaking,
a smooth quasiprojective variety
$E'_G$ carrying an action of $G$ together with
a dominant morphism $\psi\,\colon E'_G\,\longrightarrow\, X$ such
that $E'_G$ is a principal $G$-bundle over $X\setminus D$,
and the isotropy groups corresponding to the action of $G$
on the fibers of $\psi$
over the points of $D$ are finite (for a precise definition see Section
\ref{parabolic}). Let $P$ be a proper parabolic subgroup of $G$ and
$\chi\,\colon \,P\,\longrightarrow\, {\mathbb G}_m\,=\,
{\mathbb C}^\ast$ a nontrivial antidominant character.
Denote by $N$ the least common multiple of the orders of the isotropy
groups for the action of $G$ on $E'_G$. Then the character
$\chi^N$ defines a line bundle over $E'_G/P$. We show that the
parabolic bundle is semistable if and only if this line bundle
over $E'_G/P$ is numerically effective (Proposition \ref{prop.-1}).

\section{Criterion for semistability over a curve}\label{s.12}

Let $G$ be a simple linear algebraic group
defined over the field of complex
numbers. A parabolic subgroup of $G$ is a connected Zariski closed
proper subgroup $P\, \subsetneq\, G$ such that $G/P$ is a complete
variety. The quotient map $G\, \longrightarrow\, G/P$
defines a holomorphic principal $P$--bundle over $G/P$.
We recall that a character $\chi$ of $P$ is called
\textit{antidominant} if
the line bundle $(G\times {\mathbb C}_\chi)/P$ over $G/P$,
associated to this $P$--bundle for the character $\chi$,
is numerically effective; here the action of any $p\in P$
sends $(g\, ,\lambda)\, \in\, G\times{\mathbb C}$ to $(gp\, ,
\chi(p)^{-1}\lambda)$. Note that the character group of $P$
is identified with a finite index subgroup of $\text{Pic}(G/P)$.

Let $X$ be a connected smooth projective curve
defined over $\mathbb C$.
We recall from \cite[page 129, Definition 1.1]{Ra} that a
holomorphic $G$--bundle $E_G$ over the curve $X$
is called \textit{semistable} if for every reduction
of structure group of $E_G$
$$
\sigma\,\colon X\, \longrightarrow\, E_G/P
$$
to a maximal parabolic subgroup $P$ of $G$ one has
$$
\text{degree}(\sigma^\ast T_{\text{rel}})\,\ge\, 0\, ,
$$
where $T_{\text{rel}}$ is the relative tangent bundle
over $E_G/P$
for the natural projection $E_G/P\, \longrightarrow\,
X$. Alternatively, one can say that $E_G$ is semistable
if for every triple
$(P\, , E_P\, ,\chi)$, where $P\, \subset\, G$ is a parabolic
subgroup, $E_P\,\subset\, E_G$ a reduction of structure group of
$E_G$ to $P$
and $\chi$ an antidominant character of $P$, the associated
line bundle $E_P(\chi)\, =\, (E_P\times{\mathbb C})/P$ over $X$
is of nonnegative degree \cite[page 131, Lemma 2.1]{Ra}.

The following lemma will be needed in the proof of Proposition
\ref{pro.11}.

\begin{lemma}\label{lem1}
A principal $G$--bundle $E_G$ over $X$ is semistable if and
only if there is a nontrivial finite dimensional complex
$G$--module $V$ such that the vector bundle $E_G(V)\,:=\,
(E_G\times V)/G$ over $X$ associated to $E_G$ for $V$
is semistable.
\end{lemma}

\begin{proof}
If $E_G$ is semistable, from \cite[page 285, Theorem
3.18]{RR} it follows immediately that the associated vector
bundle $E_G(V)\,:=\, (E_G\times V)/V$ is semistable.

Now assume that $E_G(V)$ is semistable, where $V$ is a
nontrivial finite dimensional complex $G$--module.

Take a maximal parabolic subgroup $Q\, \subset\, G$.
There is a parabolic subgroup $P_1\, \subset\, \text{SL}(V)$ such
that $Q\, =\, P_1\bigcap G$ (so $G/Q$ is embedded in
$\text{SL}(V)/P_1$), and some positive multiple of any given
ample line bundle over $G/Q$ is the restriction
of some ample line bundle on $\text{SL}(V)/P_1$.
To prove this assertion,
consider the Jordan--H\"older filtration of $V$ for the
action of the unipotent radical of $Q$
on $V$ (the unipotent radical
acts trivially on each successive quotient
of the Jordan--H\"older filtration); since the unipotent
radical is a normal subgroup of $Q$, the action of $Q$ on $V$
preserves this filtration. The above parabolic subgroup
$P_1$ can be taken to be the subgroup of $\text{SL}(V)$ that
preserves this filtration of $V$. Note that
since $\text{Pic}(G/Q)\, \cong\, {\mathbb Z}$,
the restriction to $G/Q$ of any ample line bundle over
$\text{SL}(V)/P_1$ is a positive multiple of the (unique)
ample generator of $\text{Pic}(G/Q)$. 

Therefore, if
$\sigma\, :\, X\, \longrightarrow\, E_G/Q$ is a reduction
of structure group of $E_G$ to $Q$, then
$\text{degree}(\sigma^*T_{G/Q})$ is a positive multiple of
$\text{degree}(\sigma^*_1 T_{\text{SL}(V)/P_1})$, where
$$
\sigma_1\, :\, X\, \longrightarrow\, E_G(\text{SL}(V))/P_1
$$ 
is the reduction of structure group to
$P_1\, \subset\, \text{SL}(V)$, defined by $\sigma$,
of the principal $\text{SL}(V)$--bundle obtained
by extending the structure group of $E_G$
(by the injective homomorphism $V\,\longrightarrow\,
\text{SL}(V)$), and $T_{G/Q}$ (respectively,
$T_{\text{SL}(V)/P_1}$) is
the relative tangent bundle over $E_G/Q$ (respectively,
$E_G(\text{SL}(V))/P_1$) for the
natural projection to $X$.
Consequently, the $G$--bundle $E_G$ is semistable if
the vector bundle $E_G(V)$ is so. This completes the
proof of the lemma.
\end{proof}

For a principal $G$--bundle $E_G$ over $X$ the quotient
$E_G/P$ is a fiber bundle over $X$ with fiber $G/P$, and furthermore,
the projection $E_G\, \longrightarrow\, E_G/P$ defines a principal
$P$--bundle over $E_G/P$. For any character $\chi$ of $P$, let
\begin{equation}\label{lc}
L_\chi\, :=\, (E_G\times {\mathbb C}_\chi)/P
\end{equation}
be the line bundle over $E_G/P$ associated to this $P$--bundle
for the character $\chi$ of $P$; the action of $P$ on
$E_G\times {\mathbb C}_\chi$ is defined as before.

\begin{proposition}\label{pro.11}
A principal $G$--bundle $E_G$ over $X$ is semistable if and only
if there is a pair $(P\, ,\chi)$, where $P\, \subset\, G$ is
a proper parabolic subgroup and
$$
\chi\, :\, P\, \longrightarrow\, {\mathbb G}_m\, =\,{\mathbb C}^*
$$
a nontrivial antidominant character, such that the associated
line bundle $L_\chi$ over $E_G/P$ (see
Eq. \eqref{lc}) is numerically effective.
\end{proposition}

\begin{proof}
First assume that $E_G$ is semistable.

Let
\begin{equation}\label{df}
f\, :\, E_G/P\, \longrightarrow\, X
\end{equation}
be the natural
projection. Any fiber of $f$ is isomorphic to $G/P$, and
the restriction of $L_\chi$, defined in Eq. \eqref{lc},
to a fiber of $f$ is isomorphic to the associated line bundle
$(G\times {\mathbb C})/P$ over $G/P$ corresponding to $\chi$. As
$\chi$ is antidominant, we conclude that the
restriction of $L_\chi$ to a fiber of $f$ is
numerically effective.

Let $Y_0\, \subset\, E_G/P$ be an irreducible
curve which is not contained in any fiber of $f$. Let
$\iota \, :\, Y\, \longrightarrow\, Y_0$ be the normalization.
Therefore, $Y$ is a smooth curve, and the map
\begin{equation}\label{fY}
f_Y \, :=\, f\circ\iota\, :\, Y\, \longrightarrow\, X
\end{equation}
makes $Y$
a (possibly ramified) covering of $X$. Note that there
is a tautological section
$$
\sigma\, :\, Y\, \longrightarrow\, (f^*_Y E_G)/P\, \cong\,
f^*_Y (E_G/P)
$$
that sends any $y\in Y$ to the point of $f^*_Y (E_G/P)$ defined
by $\iota(y)\,\in\, E_G/P$.
Let $E^Y_P\, \subset\, f^*_Y E_G$ be the reduction
of structure group to $P$, of the $G$--bundle $f^*_Y E_G$ over $Y$,
defined by $\sigma$.

Note that the $P$--bundle $E^Y_P$ is identified with
the principal $P$--bundle $\iota^*E_G$,
namely the pullback to $Y$ of the principal
$P$--bundle $E_G\, \longrightarrow\, E_G/P$.
Therefore, the pullback $\iota^*L_\chi$
over $Y$ is naturally identified with the
line bundle $E^Y_P(\chi)$ associated to $E^Y_P$ for the
character $\chi$ of $P$. Consequently,
to prove that $L_\chi$ is numerically effective it is
enough to show that the $G$--bundle $f^*_Y E_G$ is semistable,
where $f_Y$ is defined in Eq. \eqref{fY}.

Any finite index subgroup of a finitely presented group $\Gamma$
contains a normal subgroup of $\Gamma$ of finite index. Therefore
we have a (possibly ramified) covering $Y'\, \longrightarrow\,
Y$ such that the composition
$$
Y'\, \longrightarrow\, Y\, \longrightarrow\, X\, ,
$$
which we will denote by $f_1$, is a finite
Galois covering. Let $\Gamma_0\, =\, \text{Gal}(Y'/X)$ denote
the Galois group for $f_1$. It is easy to see
that if $f^*_1 E_G$ is semistable, then $f^*_Y E_G$ is semistable
(the pull back to $Y'$
of a destabilizing reduction of $f^*_Y E_G$ destabilizes
$f^*_1 E_G$).

If $f^*_1 E_G$ is not semistable, then it admits a Harder--Narasimhan
reduction of structure group $E'_Q\, \subset\, f^*_1 E_G$
to a proper parabolic subgroup $Q\, \subset\, G$
\cite[page 694, Theorem 1]{AAB}.
The uniqueness of the Harder--Narasimhan reduction
implies that the natural action of the Galois group $\Gamma_0$ on
the total space of
$f^*_1 E_G$ (lifting the action of $\Gamma_0$ on $Y'$) leaves 
the submanifold
$E'_Q$ invariant. Consequently, the reduction $E'_Q
\, \subset\, f^*_1 E_G$ descends
to a reduction of structure group of $E_G$ to $Q$. Since this
reduction of structure group to $Q$ satisfies all the
conditions of a Harder--Narasimhan reduction, the $G$--bundle
$E_G$ is not semistable. In other words, we have a contradiction,
and hence $f^*_1 E_G$ is semistable. Therefore, the $G$--bundle
$f^*_Y E_G$ over $Y$ is semistable. As we saw earlier, this
implies that $L_\chi$ over $E_G/P$ is numerically effective.

To prove the converse, let $P$ and $\chi$ be such that the line
bundle $L_\chi$ over $E_G/P$ is numerically effective
and nontrivial.
Let $\theta_\chi\, :=\, (G\times {\mathbb C}_{\chi})/P$
be the associated line bundle over $G/P$, defined
by $\chi$, associated to the $P$--bundle given by
the natural projection $G\, \longrightarrow\, G/P$.
Since $\chi$ is antidominant, the line bundle $\theta_\chi$ is
numerically effective, and $\theta_\chi$ is nontrivial
as $\chi$ is so.

Let
\begin{equation}\label{def.-V}
V\, :=\, H^0(G/P,\, \theta_\chi)
\end{equation}
be the irreducible $G$--module (the irreducibility
is a part of the Borel--Weil theorem; see
\cite[page 21, Theorem 2.9.1(ii)]{PS}). Note that since
$\theta_\chi$ is nontrivial, the $G$--module $V$ is also
nontrivial. We also know that for any $n\geq 1$,
the irreducible $G$--module 
\begin{equation}\label{def.-Vn}
V_n\, :=\, H^0(G/P,\, \theta^{\otimes n}_\chi)
\end{equation}
is a direct summand of the $G$--module $V^{\otimes n}$.
To see this note that there is a
set--theoretic map from $V$ to $V_n$ that sends
any section
$$
s \,\in\, H^0(G/P,\, \theta_\chi)
$$
to $s^{\otimes n}\,\in\, H^0(G/P,\, \theta^{\otimes n}_\chi)$.
The linear span of the image of this map is a quotient of
the symmetric product $\text{Sym}^n(H^0(G/P,\,
\theta_\chi))$. Therefore, the $G$--module
$V_n$ is a direct summand of $\text{Sym}^n(V)$.
Since the $G$--module $\text{Sym}^n(V)$ is a direct
summand of $V^{\otimes n}$, we conclude
that $V_n$ is a direct summand of $V^{\otimes n}$.

Let $E_G(V)\, :=\, (E_G\times V)/G$ be the vector bundle
over $X$ associated to $E_G$ for the $G$--module $V$.
In view of Lemma \ref{lem1}, to prove that $E_G$ is
semistable it suffices to show that the vector bundle
$E_G(V)$ is semistable.

The vector bundle $E_G(V)$
is clearly identified with the direct
image $f_*L_\chi$, and similarly $E_G(V_n)$ is identified
with $f_*L^{\otimes n}_\chi$, where $V$ and $V_n$
are defined in Eq. \eqref{def.-V} and Eq. \eqref{def.-Vn}
respectively, and $f$ is the projection in Eq. \eqref{df}.
The group $G$ being simple does not have any nontrivial
character. Hence the degree of every vector bundle
associated to $E_G$ is zero.

Assume that the vector bundle $E_G(V)\, \cong\,
f_*L_\chi$ is not semistable. Let $F\, \subset\, f_*L_\chi$
be the (unique) maximal semistable subbundle $E_G(V)$. So
$F$ is the first term of the Harder--Narasimhan filtration of
$E_G(V)$. Since $\text{degree}(E_G(V))\, =\, 0$, we have
$\text{degree}(F)\, >\, 0$. For any $n \geq 1$, let
\begin{equation}\label{def.-Fn}
F_n\, \subset\, f_*L^{\otimes n}_\chi\, \cong\, E_G(V_n)
\end{equation}
be the coherent subsheaf generated by $\{v^{\otimes n}\}_{v\in F}$.
This subsheaf $F_n$ is a quotient of the
symmetric product $\text{Sym}^n (F)$. Note that
$$
\frac{\text{degree}(\text{Sym}^n (F))}{\text{rank}(\text{Sym}^n
(F))}\, =\, \frac{n\cdot \text{degree}(F)}{\text{rank}(F)}\, ,
$$
and since $F$ is semistable, the symmetric product
$\text{Sym}^n(F)$ is semistable \cite[page 285, Theorem 3.18]{RR}.
Consequently, for the quotient $F_n$ of $\text{Sym}^n (F)$ we have
\begin{equation}\label{def.-de.-ra.}
\frac{\text{degree}(F_n)}{\text{rank}(F_n)}\, \geq \, 
\frac{n\cdot \text{degree}(F)}{\text{rank}(F)}\, .
\end{equation}

Let $p_0/q_0$, with $p_0,q_0\, \in\, {\mathbb N}$, be a
positive rational number such that
$$
\frac{p_0}{q_0} < \frac{\text{degree}(F)}{\text{rank}(F)}
$$
(recall that $\frac{\text{degree}(F)}{\text{rank}(F)}\, >\, 0$).
Fix a point $x_0\,\in\, X$. For any integer
$m\geq 1$, consider the vector bundle
$$
W_m \, :=\, F_{mq_0}\otimes {\mathcal O}_X(-mp_0 x_0)
$$
over $X$
where $F_{mq_0}$ is defined in Eq. \eqref{def.-Fn}. Using 
Eq. \eqref{def.-de.-ra.} we have
$$
\frac{\text{degree}(W_m)}{\text{rank}(W_m)}\, \geq\,
mq_0(\frac{\text{degree}(F)}{\text{rank}(F)}- \frac{p_0}{q_0})
\, >\, 0\, .
$$
In view of this inequality, if $m$ is such that
$mq_0(\frac{\text{degree}(F)}{\text{rank}(F)}- \frac{p_0}{q_0})
\,>\, \text{genus}(X) -1$, then the Riemann--Roch theorem gives
$$
\dim H^0(X,\, W_m)\, \geq\,\chi(W_m)\, \geq\,
\text{rank}(W_m)(mq_0(\frac{\text{degree}(F)}{\text{rank}(F)}-
\frac{p_0}{q_0})-\text{genus}(X)+1) \, >\, 0\, .
$$

Take an integer $m$ such that
$mq_0(\frac{\text{degree}(F)}{\text{rank}(F)}- \frac{p_0}{q_0})
\,>\, \text{genus}(X) -1$, and take a nonzero section
\begin{equation}\label{s0}
0\,\not=\, s_0\,\in\, H^0(X,\, W_m)\,\subset\,
H^0(X,\, E_G(V_{mq_0})\otimes {\mathcal O}_X(-mp_0 x_0))\, .
\end{equation}
Using the isomorphism
$f_*L^{\otimes n}_\chi\, \cong\, E_G(V_n)$ and the
projection formula we have
$$
H^0(X,\, E_G(V_{mq_0})\otimes {\mathcal O}_X(-mp_0 x_0))
\,=\, H^0(E_G/P,\, L^{\otimes mq_0}_\chi\otimes
f^*{\mathcal O}_X(-mp_0 x_0))\, .
$$
Let
\begin{equation}\label{sp0}
s'_0 \, \in\, H^0(E_G/P,\, L^{\otimes mq_0}_\chi\otimes
f^*{\mathcal O}_X(-mp_0 x_0))
\end{equation}
be the section defined by $s_0$ (in Eq. \eqref{s0})
using the above isomorphism. Set
$$
\text{Div}(s'_0) \, \subset\, E_G/P
$$
to the effective divisor defined by the above section $s'_0$.

Let $d_0$ be the largest integer such that the cup product
$c_1(\theta_\chi)^{d_0} \, \in\, H^{2d_0}(G/P,\,
{\mathbb Q})$ is nonzero. Since $\theta_\chi$ is
numerically effective but nontrivial, there is
a parabolic subgroup $P'\, \supseteq\, P$ of $G$ such that
$\theta_\chi$ is the pull back of an ample line bundle over
$G/P'$ using the natural projection of $G/P$ to $G/P'$.
The above defined integer $d_0$ is the dimension of $G/P'$.

Let $\chi'$ be the character of $P'$ that restricts to
$\chi$ on $P$. Equivalently, the line bundle $\theta_{\chi'}$
over $G/P'$ defined by $\chi'$ satisfies the condition that
$\theta_{\chi'}$ pulls back to $\theta_{\chi}$.

Let
$$
\phi\, :\, E_G/P\, \longrightarrow\, E_G/P'
$$
be the natural projection and $L_{\chi'}\,:=\,
(E_G\times{\mathbb C}_{\chi'})/P'$ the line bundle over
$E_G/P'$ defined by $\chi'$. So, we have $\phi^*L_{\chi'}\,
\cong \, L_\chi$.

Since $L_\chi$ is numerically effective and $\phi$ is
surjective, we conclude that $L_{\chi'}$
is numerically effective \cite[page 360, Proposition 2.3]{Fu}.

Since $\phi_{*}{\mathcal O}_{E_G/P} \, \cong\,
{\mathcal O}_{E_G/P'}$, we have
$$
\phi_*(L^{\otimes mq_0}_\chi\otimes
f^*{\mathcal O}_X(-mp_0 x_0))
\, \cong\, L^{\otimes mq_0}_{\chi'}\otimes
h^*{\mathcal O}_X(-mp_0 x_0)\, ,
$$
where
\begin{equation}\label{de.h}
h\, :\, E_G/P'\, \longrightarrow\, X
\end{equation}
is the natural projection. Therefore,
$$
H^0(E_G/P,\, L^{\otimes mq_0}_\chi\otimes
f^*{\mathcal O}_X(-mp_0 x_0))\, \cong\,
H^0(E_G/P',\, L^{\otimes mq_0}_{\chi'}\otimes
h^*{\mathcal O}_X(-mp_0 x_0))\, .
$$
Consequently, the section $s'_0$ of $L^{\otimes mq_0}_\chi
\otimes f^*{\mathcal O}_X(-mp_0 x_0)$ constructed in
Eq. \eqref{sp0} corresponds to a section
\begin{equation}\label{def.-z-0}
\zeta_0 \, \in\, H^0(E_G/P',\, L^{\otimes mq_0}_{\chi'}
\otimes h^*{\mathcal O}_X(-mp_0 x_0))\, ,
\end{equation}
and the divisor $\text{Div}(s'_0)$ on $E_G/P$ coincides
with the divisor $\phi^{-1}(\text{Div}(\zeta_0))$.

We will show that $(c_1(L_{\chi'}))^{d_0+1}
\, \in\, H^{2d_0+2}(E_G/P',
\, {\mathbb Z})\, \cong\, {\mathbb Z}$ vanishes.

The topological isomorphism classes of $G$--bundles
over a compact connected Riemann surface
are parametrized by the fundamental group
$\pi_1(G)$ \cite[page 142, Proposition 5.1]{Ra}.
Since $G$ is simple, we know that
$\pi_1(G)$ is a finite group.
Therefore, if $\beta\, :\, X'\, \longrightarrow\, X$
is a finite cover (possibly ramified) of degree
$\# \pi_1(G)$ with $X'$ connected,
then $\beta^*E_G$ is topologically trivial;
here $\# \pi_1(G)$
is the cardinality of $\pi_1(G)$. Fix a covering
$$
\beta\, :\, X'\, \longrightarrow\, X
$$
such that $\beta^*E_G$ is topologically trivial.

Let $\widehat{\beta} \, :\, \beta^*E_G/P'\, \longrightarrow\,
E_G/P'$ be the natural projection over $\beta$. So
$\widehat{\beta}^*L_{\chi'}$ is isomorphic to the line bundle
$\beta^*E_G(\chi')\, :=\,
(\beta^*E_G\times {\mathbb C}_{\chi'})/P'$
(the line bundle over $\beta^*E_G/P'$ associated to the
principal $P'$--bundle defined by $\beta^*E_G\,
\longrightarrow\, \beta^*E_G/P'$ for the character
$\chi'$). If we fix a topological isomorphism
$$
\tau\, :\, \beta^*E_G\, \longrightarrow\, X'\times G
$$
of $\beta^*E_G$ with the trivial $G$--bundle over $X'$, then the
line bundle $\beta^*E_G(\chi')$ over $\beta^*E_G/P'$ is
identified with the pullback $(p_{G/P'}\circ\tau)^*
{\mathcal L}_{P'}$, where $p_{G/P'}$ is the composition of the
projection of $X'\times G\, \longrightarrow\, G$ with the
projection $G\, \longrightarrow\, G/P$, and ${\mathcal L}_{P'}$
is the line bundle over $G/P'$ defined by $\chi'$.
Therefore, we have
$$
(c_1(\beta^*E_G(\chi')))^{d_0+1}\, =\, 0\, .
$$
Since $(c_1(\beta^*E_G(\chi')))^{d_0+1}\, \in\, {\mathbb Q}$
coincides with $\# \pi_1(G)\cdot (c_1(L_{\chi'}))^{d_0+1}$
(the degree of $\widehat{\beta}$ coincides with the degree
of $\beta$, and the degree of $\beta$ is $\# \pi_1(G)$),
this implies that $(c_1(L_{\chi'}))^{d_0+1}\, =\, 0$.

Since the cohomology class
$$
c_1(L^{\otimes mq_0}_{\chi'}))^{d_0+1} \,=\,
(mq_0)^{d_0+1}(c_1(L_{\chi'}))^{d_0+1}\, \in\,
H^{2d_0+2}(E_G/P', \, {\mathbb Z})
$$
is zero, it can be shown that the integer
$$
(c_1(L^{\otimes mq_0}_{\chi'}))^{d_0}
c_1(L^{\otimes mq_0}_{\chi'}\otimes
h^*{\mathcal O}_X(-mp_0 x_0))\, \in\,
H^{2d_0+2}(E_G/P', \, {\mathbb Z})\, \cong\, {\mathbb Z}
$$
coincides with
$$
-mp_0 (c_1(\theta_{\chi'}))^{d_0}\, \in\, H^{2d_0}(G/P', \,
{\mathbb Z}) \, \cong\, {\mathbb Z}\, ,
$$
where $d_0$, as before, is the dimension of $G/P'$ and
$h$ is the projection in Eq. \eqref{de.h}. Indeed, we have
$c_1(L^{\otimes mq_0}_{\chi'}\otimes
h^*{\mathcal O}_X(-mp_0 x_0))\, =\,
c_1(L^{\otimes mq_0}_{\chi'})
-mp_0 [h^{-1}(x_0)]$, where $[h^{-1}(x_0)]$ is cohomology
class defined by $h^{-1}(x_0)$ using Poincar\'e duality.
Hence
$$
(c_1(L^{\otimes mq_0}_{\chi'}))^{d_0}
c_1(L^{\otimes mq_0}_{\chi'}\otimes
h^*{\mathcal O}_X(-mp_0 x_0))
$$
$$
=\, c_1(L^{\otimes mq_0}_{\chi'}))^{d_0+1} -
mp_0 (c_1(\theta_{\chi'}))^{d_0}\, =\,
-mp_0 (c_1(\theta_{\chi'}))^{d_0}\, .
$$

Since $\theta_{\chi'}$ is ample, we have
$(c_1(\theta_{\chi'}))^{d_0}\, >\, 0$. Therefore,
$$
(c_1(L^{\otimes mq_0}_{\chi'}))^{d_0}
c_1(L^{\otimes mq_0}_{\chi'}\otimes
h^*{\mathcal O}_X(-mp_0 x_0))\, <\, 0\, .
$$
This means that the restriction of the line bundle
$L^{\otimes mq_0}_{\chi'}$ to the effective divisor
$\text{Div}(\zeta_0)$ (constructed in Eq. \eqref{def.-z-0})
on $E_G/P'$ is of negative degree
(recall that $\text{Div}(\zeta_0)$ defines $L^{\otimes
mq_0}_{\chi'}\otimes h^*{\mathcal O}_X(-mp_0 x_0)$); by degree
we mean the top exterior product of the first Chern class. But
this contradicts the fact that the line bundle
$L^{\otimes mq_0}_{\chi'}$ is numerically effective.

Therefore, $E_G$ is semistable, and the proof of the proposition
is complete.
\end{proof}

\section{Principal bundles over projective manifolds}\label{s.13}

Let $M$ be a connected complex projective manifold
of complex dimension $d$. Fix a very ample line bundle
$\xi \, \in\, \text{Pic}(M)$ on $M$. For
a coherent sheaf $F$ on $M$, define the \textit{degree}
$$
\text{degree}(F)\, :=\, (c_1(F)\cup c_1(\xi)^{d-1})\cap [M]
\,\in\, {\mathbb Z}\, .
$$
If $F$ is a holomorphic vector bundle defined over a
Zariski open dense subset $U\, \subset\, M$ such
that the complement $M\setminus U$ is
of (complex) codimension at least two, then define
$$
\text{degree}(F)\, :=\, \text{degree}(\iota_* F)\, ,
$$
where $\iota\, :\, U\, \hookrightarrow\, M$ is the
inclusion map; note that $\iota_* F$ is a coherent sheaf on $M$.

Let $G$ be a connected reductive linear algebraic group defined
over $\mathbb C$. Let $Z(G)\, \subset\, G$ be the center of $G$.
A principal $G$--bundle $E_G$ over $M$ is called \textit{semistable}
if for any holomorphic
reduction of structure group $E_P\, \subset\, E_G\vert_U$ to any
parabolic subgroup $P\, \subset\, G$ over
some Zariski open subset $U\, \subset M$,
with $\text{codim}_{\mathbb C}(M\setminus U) \, \geq \,
2$, and for any nontrivial antidominant
character $\chi$ of $P$ which is trivial of $Z(G)$, the
associated line bundle $E_P(\chi) \, :=\,
(E_P\times{\mathbb C})/P$ over $U$ satisfies the condition
$$
\text{degree}(E_P(\chi))\, \geq \, 0
$$
(see \cite{Ra}, \cite{AB}); in the above definition of
$E_P(\chi)$, the action of $g\, \in\, P$ sends any $(z\, ,
c)\,\in\,E_P\times{\mathbb C}$ to $(zg\, ,\chi(g^{-1})c)$.

Equivalently, $E_G$ is semistable if for every
reduction $\sigma\,\colon\, M\, \longrightarrow
\, (E_G/P)\vert_U$ to a (proper) maximal parabolic subgroup
$P$ of $G$, over any Zariski open subset $U\, \subset\, M$
of the above type, one has
$$
\text{degree}(\sigma^\ast T_{\text{rel}})\,\ge\, 0\, ,
$$
where $T_{\text{rel}}$ is the relative tangent bundle of the
projection $E_G/P\, \longrightarrow\, M$ (see \cite[page 131,
Lemma 2.1]{Ra}).

Note that in the special case where $\dim M\,=\, 1$
and $G$ is a simple group the above definition
of semistability coincides with the one given in Section
\ref{s.12}.

Let $\mathfrak g$ be the Lie algebra of $G$. For a holomorphic
$G$--bundle $E_G$ over $M$, let
$$
\text{ad}(E_G)\, :=\, (E_G\times {\mathfrak g})/G
$$
be the adjoint vector bundle over $M$; the action of
$g\, \in\, G$ sends any $(z\, ,v)\, \in\,
E_G\times{\mathfrak g}$ to $(zg\, ,\text{ad}(g^{-1})(v))$.
It is known that $E_G$ is semistable if and only if the
vector bundle $\text{ad}(E_G)$ is semistable
\cite[page 214, Proposition 2.10]{AB}.

Let $P\, \subset\, G$ a parabolic subgroup and $\chi$
a character of $P$. So we have an associated line bundle
\begin{equation}\label{l.-chi}
L_\chi \, :=\, (E_G\times {\mathbb C}_\chi)/P
\end{equation}
over $E_G/P$; the action of $g\, \in\, P$ sends any $(z\, ,c)\,
\in\, E_G\times{\mathbb C}$ to $(zg\, ,\chi(g^{-1})c)$. Let
\begin{equation}\label{de.-pr.}
\phi\, :\, E_G/P\, \longrightarrow\, M
\end{equation}
be the natural projection. Any fiber of $\phi$ is isomorphic
to $G/P$.

\begin{remark}\label{rem.-0}
{\rm Let $G$ be a simple linear algebraic group
over $\mathbb C$. If $\dim M\, =\, 1$,
then a $G$--bundle $E_G$ over $M$ is semistable if and
only if the line bundle $L_\chi$ in Eq. \eqref{l.-chi}
is numerically effective for some parabolic subgroup
$P\, \subsetneq\, G$ and some nontrivial antidominant character
$\chi$ of $P$ (Proposition \ref{pro.11}). On the other hand,
the main theorem of \cite{MR} says that a vector bundle $V$ over
$M$ is semistable if and only if the restriction $V\vert_C$ is
semistable, where $C\, \subset\, M$ is a general complete
intersection curve of sufficiently large degree hypersurfaces in
$M$ (corresponding to $\xi$); see \cite[page 221, Theorem 6.1]{MR}.
In particular, $\text{ad}(E_G)$ is semistable if and only if
$\text{ad}(E_G)\vert_C$ is semistable for any such $C$. Therefore,
using the criterion (in Proposition \ref{pro.11}) for semistability
of a $G$--bundle over a Riemann surface it follow immediately that
$E_G$ is semistable if and only if the line bundle
$L_\chi\vert_{\phi^{-1}(C)}$ is numerically effective, where
$L_\chi$ (respectively, $\phi$) is defined in Eq. \eqref{l.-chi}
(respectively, Eq. \eqref{de.-pr.}) and $C\, \subset\, M$ is a
general complete intersection curve of sufficiently
large degree hypersurfaces in $M$ (corresponding to $\xi$).}
\end{remark}

\begin{theorem}\label{thm.h.d.}
Let $E_G$ be a holomorphic principal $G$--bundle over
a connected projective manifold $M$, where $G$ is a
simple linear algebraic group defined over $\mathbb C$.
Fix a proper parabolic subgroup $P\, \subset\, G$ and a
nontrivial antidominant character $\chi$ of $P$.
The following two statements are equivalent:
\begin{enumerate}

\item{} The associated line bundle $L_\chi\, :=\,
(E_G\times {\mathbb C}_\chi)/P$ over $E_G/P$ defined by $\chi$
is numerically effective.

\item{} The $G$--bundle $E_G$ is semistable and the second Chern
class
$$
c_2({\rm ad}(E_G)) \, \in\, H^4(M,\, {\mathbb Q})
$$
vanishes.

\end{enumerate}
\end{theorem}

\begin{proof}
First assume that the line bundle
$L_\chi$ is numerically effective. Let $C\, \subset\,
M$ be a connected smooth complex projective curve. Since $L_\chi$
is numerically effective, the restriction of $L_\chi$
to $\phi^{-1}(C)$ is numerically effective, where $\phi$
is the projection in Eq. \eqref{de.-pr.}.
So using Proposition \ref{pro.11} it follows that
the restriction $E_G\vert_C$ of $E_G$ to $C$ is semistable.

Since the restriction of $E_G$ to every connected smooth complex
projective curve in the variety $M$
is semistable, we conclude that the restriction of
$\text{ad}(E_G)$ to every connected
smooth curve in $M$ is semistable;
recall from Remark \ref{rem.-0} that
a $G$--bundle is semistable if and only if
its adjoint bundle is semistable (see \cite[page 214, Proposition
2.10]{AB} for proof). From this it follows
that the vector bundle $\text{ad}(E_G)$ is semistable. Indeed,
if $F\, \subset\, \text{ad}(E_G)$ is a subsheaf contradicting
the semistability condition, then the restriction to $F$
to a smooth curve $C\, \subset\, M$ satisfying the two conditions
\begin{enumerate}
\item{} $C$ is contained in the Zariski open
dense subset over which
$F$ is a subbundle of $\text{ad}(E_G)$ (the complement of this
open set is of codimension at least two),
\item{} $C$ is a complete intersection of hypersurfaces
in $M$ for the given very ample line bundle $\xi$.
\end{enumerate}
contradicts the semistability condition of $\text{ad}(E_G)\vert_C$.

As $\text{ad}(E_G)$ is semistable, we conclude that $E_G$ is
semistable.
We still need to show that $c_2(\text{ad}(E_G))\, =\, 0$ to be
able to conclude that statement (1) implies statement (2).

For a vector bundle $E$ over $M$, let ${\mathbb P}(E)$ denote
the projective bundle over $M$ defined by the one--dimensional
quotients of the fibers of $E$. The tautological line bundle
over ${\mathbb P}(E)$ will be denoted by ${\mathcal O}_{{\mathbb
P}(E)}(1)$.
A vector bundle $E$ over $M$ is called \textit{numerically effective}
if the line bundle ${\mathcal O}_{{\mathbb P}(E)}(1)$ over ${\mathbb
P}(E)$ is numerically effective.

We will use properties of
numerically effective vector bundles proved in \cite{DPS}.
In \cite{DPS}, a holomorphic vector bundle $E$ over a compact
complex Hermitian manifold $X$ is called \textit{numerically
effective} if ${\mathcal O}_{{\mathbb P}(E)}(1)$ over
${\mathbb P}(E)$ admits Hermitian connections
with arbitrary small negative part of curvature
(see \cite[page 297, Definition 1.2]{DPS}). For
projective manifolds, the above two definitions of
numerically effectiveness coincide.
Indeed, if $E$ is numerically effective in the sense of
\cite[Definition 1.2]{DPS}, then it is clearly
numerically effective when defined using curves in
${\mathbb P}(E)$; see the comment in \cite[page 297]{DPS}
following Definition 1.2. For the converse direction,
take any connected complex projective manifold $X$, and
fix an ample line bundle $L$ on $X$. Also, fix
a Hermitian connection on $L$ with positive curvature. Let
$\Omega_L$ be the curvature of this
Hermitian connection on $L$, which is
positive by assumption. If $\eta$ is a numerically
effective line bundle on $X$, then
$\eta^{\otimes n}\otimes L$ is ample for each $n\geq 1$
(numerically effectiveness is defined using curves). Take
a Hermitian connection on $\eta^{\otimes n}\otimes L$ with positive
curvature. This connection on $\eta^{\otimes n}\otimes L$ and the
given connection on $L$ together define a Hermitian
connection on $\eta$. Indeed, as $\eta^{\otimes
n}\, =\, (\eta^{\otimes n}\otimes L)\otimes L^*$, there
is an induced connection on $\eta^{\otimes n}$, which in turn
induces a connection on $\eta$. The
curvature of this connection on $\eta$ is bounded below by
$-\Omega_L/n$ (recall that the curvature of the connection on
$\eta^{\otimes n}\otimes L$ is positive). Therefore, the line
bundle $L$ is numerically effective in the sense of \cite{DPS}.

Let
\begin{equation}\label{def.-f-ad}
f\, :\, C\, \longrightarrow\, M
\end{equation}
be a holomorphic map from a connected smooth complex projective
curve. We will show that the vector bundle
$f^*\text{ad}(E_G)$ over $C$ is numerically effective. For this
we will first prove that the $G$--bundle $f^* E_G$ is semistable.

There is a natural map
$$
f_G\, :\, (f^* E_G)/P\, \longrightarrow\, E_G/P
$$
defined using the natural identification $(f^* E_G)/P
\,=\, f^* (E_G/P)$. The line bundle $f^*_G L_\chi$
over $(f^* E_G)/P$ is identified with the line
bundle $L'_\chi$ associated, for the character $\chi$ of $P$, to
the principal $P$--bundle over $(f^* E_G)/P$ defined by
the projection $f^* E_G\, \longrightarrow\, (f^* E_G)/P$.
Since $L_\chi$ is numerically effective, we conclude
that $f^*_G L_\chi\, =\, L'_\chi$ over $(f^* E_G)/P$ is
also numerically effective
\cite[page 360, Proposition 2.2]{Fu}. Consequently, the $G$--bundle
$f^* E_G$ is semistable (Proposition \ref{pro.11}).

To show that $f^*\text{ad}(E_G)$ over $C$ is numerically effective,
first note that the vector bundle $f^*\text{ad}(E_G)$ is semistable
as the $G$--bundle $f^* E_G$ is semistable. Now,
since $(f^*\text{ad}(E_G))^*\, \cong\, f^*\text{ad}(E_G)$
(any $G$--invariant nondegenerate symmetric bilinear form
on $\mathfrak g$ gives a nondegenerate symmetric bilinear form
on $\text{ad}(E_G)$), we conclude that
the line bundle $\bigwedge^{\text{top}}f^*\text{ad}(E_G)$ is
trivial. This implies that the line bundle
\begin{equation}\label{de.-l.b.}
{\mathcal O}_{{\mathbb P}(f^*\text{ad}(E_G))}(\dim_{\mathbb C} G)
\,:=\,{\mathcal O}_{{\mathbb P}(f^*\text{ad}(E_G))}(1)^{\otimes
(\dim_{\mathbb C} G)}
\end{equation}
over the total space of the projective bundle
${\mathbb P}(f^*\text{ad}(E_G))\, \longrightarrow\, C$
is identified with the top exterior product of the relative
tangent bundle for the natural projection of
${\mathbb P}(f^*\text{ad}(E_G))$ to $C$. Using this and
the fact that $f^*\text{ad}(E_G)$ is semistable it
can be deduced from Proposition \ref{pro.11} that the line bundle
in Eq. \eqref{de.-l.b.} is numerically effective. To deduce this
assertion from Proposition \ref{pro.11}, consider the principal
$\text{PGL}(\mathfrak g)$--bundle
$E_{\text{PGL}(\mathfrak g)}$ over $C$ defined by the projective
bundle ${\mathbb P}(f^*\text{ad}(E_G))$; take the maximal parabolic
subgroup $Q'$ of $\text{PGL}(\mathfrak g)$ that fixes a hyperplane 
in $\mathfrak g$, and note that the anticanonical line
bundle over $\text{PGL}(\mathfrak g)/Q'$ corresponds to an
antidominant character of $Q'$. As the vector bundle
$f^*\text{ad}(E_G)$ is
semistable, Proposition \ref{pro.11} says that the relative
anticanonical bundle for the projection
${\mathbb P}(f^*\text{ad}(E_G))\, \longrightarrow\, C$
is numerically effective. Thus the line bundle in
Eq. \eqref{de.-l.b.} is numerically effective.

The line bundle in Eq. \eqref{de.-l.b.} being
numerically effective we conclude that the line bundle
${\mathcal O}_{{\mathbb P}(f^*\text{ad}(E_G))}(1)$
is numerically effective. In other words, the vector bundle
$f^*\text{ad}(E_G)$ is numerically effective.

Let
\begin{equation}\label{eq.-f0}
f_0\, :\, C\, \longrightarrow\, {\mathbb P}(\text{ad}(E_G))
\end{equation}
be a holomorphic map from a connected smooth complex projective
curve. Let
$$
p\, :\, {\mathbb P}(\text{ad}(E_G)) \, \longrightarrow\, M
$$
be the natural projection. Set $f$ in Eq. \eqref{def.-f-ad} to
be $p\circ f_0$. Note that there is a natural map from
the projective bundle over $C$
$$
f_1\, :\, {\mathbb P}(f^*\text{ad}(E_G))\,=\,
f^*{\mathbb P}(\text{ad}(E_G)) \,
\longrightarrow\, {\mathbb P}(\text{ad}(E_G))
$$
which projects to $f$. Also, there is a canonical section
$$
\sigma_1\, :\, C\, \longrightarrow\,
{\mathbb P}(f^*\text{ad}(E_G))
$$
(of the projection ${\mathbb P}(f^*\text{ad}(E_G))\,
\longrightarrow\, C$); the pair $(f_1\, , \sigma_1)$ satisfy
the following two conditions
\begin{enumerate}
\item{} $f_1\circ\sigma_1 \, =\, f_0$, where $f_0$
is defined in Eq. \eqref{eq.-f0}, and
\item{} $f^*_1{\mathcal O}_{{\mathbb P}(\text{ad}(E_G))}(1)\, =\,
{\mathcal O}_{{\mathbb P}(f^*\text{ad}(E_G))}(1)$.
\end{enumerate}
We have proved above that $f^*\text{ad}(E_G)$ is numerically
effective. Using this together with
the above two properties of $f_1$ and
$\sigma_1$ we conclude that
$$
\text{degree}(f^*_0{\mathcal O}_{{\mathbb P}(\text{ad}(E_G))}(1))
=\text{degree}((f_1\circ\sigma_1)^*{\mathcal O}_{{\mathbb
P}(\text{ad}(E_G))}(1)) = \sigma^*_1
\text{degree}({\mathcal O}_{{\mathbb P}(f^*\text{ad}(E_G))}(1))
\geq 0\, .
$$

Consequently, the vector
bundle $\text{ad}(E_G)$ over $M$ is numerically effective.

A vector bundle $E$ is called \textit{numerically flat} if
both $E$ and $E^*$ are numerically effective
\cite[page 311, Definition 1.17]{DPS}.
Since $\text{ad}(E_G)^*\, \cong\, \text{ad}(E_G)$ and
$\text{ad}(E_G)$ is numerically effective, we conclude that
$\text{ad}(E_G)$ is numerically flat. All the Chern classes
of positive degree of a numerically flat vector bundle vanish
\cite[page 311, Corollary 1.19]{DPS}; in particular, we have
$c_2(\text{ad}(E_G))\, =\, 0$.

Therefore, statement (1) in the theorem implies statement (2).
To prove the converse, assume that $E_G$ is semistable and
$c_2(\text{ad}(E_G))\, =\, 0$.

The semistability of $E_G$ implies that the vector bundle
$\text{ad}(E_G)$ is semistable
\cite[Proposition 2.10]{AB}. Since $\text{ad}(E_G)^*\, \cong\,
\text{ad}(E_G)$, we have $c_1(\text{ad}(E_G))\, =\, 0$, and it
is given that
$c_2(\text{ad}(E_G))\, =\, 0$. Since $\text{ad}(E_G)$ is semistable
with $c_1(\text{ad}(E_G))\, =\, 0\, =\, c_2(\text{ad}(E_G))$, it
admits a filtration of holomorphic subbundles
\begin{equation}\label{filt.-simpson}
0\, =\, F_0\, \subset\, F_1\, \subset\, F_2\, \subset\,\cdots
\, \subset\, F_{k-1} \, \subset\, F_k \, =\, \text{ad}(E_G)
\end{equation}
such that each $F_i/F_{i-1}$, $i\,\in\, [1\, ,k]$,
is a stable vector bundle and
$c_j(F_i/F_{i-1})\, =\, 0$ for all $j\geq 1$
and $1\leq i\leq k$ \cite[page 39, Theorem 2]{Si}. To deduce
this from \cite[page 39, Theorem 2]{Si} simply set the Higgs
field in \cite[Theorem 2]{Si} to be zero.

Since $F_i/F_{i-1}$ is stable and $c_j(F_i/F_{i-1})\, =\, 0$ for
all $j\geq 1$, a theorem due to Donaldson says that each
$F_i/F_{i-1}$ in Eq. \eqref{filt.-simpson} admits a unitary flat
connection \cite[page 231, Proposition 1]{Do}.

Let
\begin{equation}\label{def.-ne.-f}
f\, :\, C\, \longrightarrow\, E_G/P
\end{equation}
be a morphism from a connected smooth complex projective curve.
We will show that the $G$--bundle $(\phi\circ f)^*E_G$ over $C$ is
semistable, where the map $\phi$ is defined in Eq. \eqref{de.-pr.}.

To prove that $(\phi\circ f)^*E_G$ is semistable, first note that
the adjoint bundle $\text{ad}((\phi\circ f)^*E_G)$ has a filtration
$$
0\, =\, F'_0\, \subset\, F'_1\, \subset\, F'_2\, \subset\,\cdots
\, \subset\, F'_{k-1} \, \subset\, F'_k \, =\,
(\phi\circ f)^*\text{ad}(E_G)\, ,
$$
where $F'_i \, :=\, (\phi\circ f)^*F_i$ with $F_i$ is as in
Eq. \eqref{filt.-simpson}. It was noted earlier that each
$F_i/F_{i-1}$ admits a unitary flat connection. A unitary flat
connection on $F_i/F_{i-1}$ pulls back to
induce a unitary flat connection on the pullback
$(\phi\circ f)^*(F'_i/F'_{i-1})\,=\, F'_i/F'_{i-1}$. Therefore,
$F'_i/F'_{i-1}$, $i\,\in\, [1\, ,k]$,
is a polystable vector bundle of degree zero
\cite[page 560, Theorem 2(A)]{NS}. Since
$(\phi\circ f)^*\text{ad}(E_G)$ is filtered by subbundles
with each successive quotient polystable of degree zero, it
follows immediately that the vector bundle
$(\phi\circ f)^*\text{ad}(E_G)$ is semistable.
Therefore, the $G$--bundle $(\phi\circ f)^*E_G$
over $C$ is semistable.

To show that the associated line bundle $L_\chi$
over $E_G/P$ defined by $\chi$ is numerically effective,
take any map $f$ as in Eq. \eqref{def.-ne.-f}. Note
that $(\phi\circ f)^*E_G$ (the map $\phi$ is defined in
Eq. \eqref{de.-pr.})
has a natural reduction of structure group to $P
\,\subset\, G$ defined by the section
$$
\sigma'\, :\, C\, \longrightarrow\, (\phi\circ f)^*(E_G/P)\,
\cong\, ((\phi\circ f)^*E_G)/P
$$
that sends any $c\in C$ to $f(c)\, \in\,
(E_G/P)_{\phi\circ f(c)} \,=\, (\phi\circ f)^*(E_G/P)_c$.
There is a natural map
$$
f_1\, :\, (\phi\circ f)^*(E_G/P)\, \longrightarrow\, E_G/P
$$
such that $f_1\circ\sigma' \, =\, f$ and the line bundle
$f^*_1 L_\chi$ over $(\phi\circ f)^*(E_G/P)$ is identified with
$((\phi\circ f)^*E_G)(\chi)$, the
line bundle over $(\phi\circ f)^*(E_G/P)$
associated, for the character $\chi$ of $P$,
to the principal $P$--bundle over $((\phi\circ f)^*E_G)/P$
defined by the natural projection $(\phi\circ f)^*E_G\,
\longrightarrow\, ((\phi\circ f)^*E_G)/P$.

Now, since $(\phi\circ f)^*E_G$ is semistable (this
was proved earlier), from
Proposition \ref{pro.11} we conclude that the line bundle
$((\phi\circ f)^*E_G)(\chi)$ over $(\phi\circ f)^*(E_G/P)$ is
numerically effective. Consequently, we have
$$
\text{degree}(f^* L_\chi)\, =\, (\sigma')^*
\text{degree}(((\phi\circ f)^*E_G)(\chi))\, \geq\, 0\, .
$$
Therefore, $L_\chi$ is numerically effective.
This completes the proof of the theorem.
\end{proof}

\begin{remark}\label{rem.-2}
{\rm Let $\kappa\, \in\, \text{Sym}^2 {\mathfrak g}^*$ be the
Killing form on the Lie algebra $\mathfrak g$ of $G$. For
a principal $G$--bundle $E_G$ over $M$, the form $\kappa$
defines a characteristic class
$$
C_\kappa(E_G) \, \in\, H^4(M,\, {\mathbb R})\, .
$$
This characteristic class can be defined as follows. For any
$C^\infty$ connection $\nabla$ on $E_G$, consider the smooth
$4$--form $\kappa (\Omega_\nabla)$ on $M$, where $\Omega_\nabla$
is the curvature of $\nabla$.
This differential form is closed and the de Rham cohomology class
defined by it does not depend on the choice of the connection
$\nabla$ on $E_G$. This cohomology class defined by the
differential form $\kappa (\Omega_\nabla)$ coincides with
$C_\kappa(E_G)$. The second Chern class
$c_2(\text{ad}(E_G))\, \in\, H^4(M,\, {\mathbb Q})$ is a
positive multiple of $C_\kappa(E_G)$
(the multiplication factor is $2m(G)$, where $m(G)$ is the
dual Coxeter number). Therefore, the condition
that $c_2(\text{ad}(E_G))\, =\, 0$
is equivalent to the condition that $C_\kappa(E_G)\,=\, 0$
(see statement (2) in Theorem \ref{thm.h.d.}).}
\end{remark}

\section{Generalization for reductive groups}\label{sec.-4}

The following proposition can be viewed as a
reformulation of Theorem \ref{thm.h.d.}.

\begin{proposition}\label{pr.-1}
Fix a pair $(P\, ,\chi)$, where $P\, \subset\, G$ is a
proper parabolic subgroup of the simple group $G$ and
$$
\chi\, :\, P\, \longrightarrow \, {\mathbb G}_m
\,=\, {\mathbb C}^*
$$
a nontrivial antidominant character. A principal $G$--bundle
$E_G$ over a projective manifold
$M$ is semistable with $c_2({\rm ad}(E_G))\, =\, 0$
if and only if for every pair of the form $(Y\, ,\psi)$,
where $Y$ is a compact connected Riemann surface and
$$
\psi\,:\, Y\,\longrightarrow\, M
$$
a holomorphic map, and every reduction
$E_P\, \subset\, \psi^*E_G$ of structure group to $P$ of the
principal
$G$--bundle $\psi^*E_G$ over $Y$, the associated line bundle
$E_P(\chi)\, =\, (E_P\times {\mathbb C}_\chi)/P$ over $Y$ is of
nonnegative degree.
\end{proposition}

\begin{proof}
If $E_G$ is semistable with $c_2(\text{ad}(E_G))\, =\, 0$,
then using a result of Simpson we saw in the proof of Theorem
\ref{thm.h.d.} that $\text{ad}(E_G)$ admits a filtration
of subbundles as in Eq. \eqref{filt.-simpson} such that each
successive quotient admits a unitary flat connection.
This filtration induces a
filtration of subbundles of the pullback $\psi^*\text{ad}(E_G)$
with the property that each successive
quotient admits a unitary flat connection (the map $\psi$ is as
in the statement of the proposition). Consequently,
$\psi^*\text{ad}(E_G)$, and hence $\psi^*E_G$, is semistable.
The semistability of $\psi^*E_G$ immediately implies that the
line bundle $E_P(\chi)$ in the statement of the proposition
is of nonnegative degree.

For the converse direction, assume that $E_P(\chi)$ is of
nonnegative degree for each reduction of structure group $E_P
\, \subset\, \psi^*E_G$ and each map $\psi$ of the above type.

Take any holomorphic map
$$
\psi_1\,:\, Y\,\longrightarrow\, E_G/P\, ,
$$
where $Y$ is compact connected Riemann surface. Set
$$
\psi\, =\, \phi\circ\psi_1\, ,
$$
where $\phi$ is the projection in Eq. \eqref{de.-pr.}.
The $G$--bundle $\psi^*E_G$ over $Y$ has a
tautological reduction of structure group to $P$
$$
\sigma\,:\, Y\, \longrightarrow\, (\psi^*E_G)/P
$$
that sends any point $y\in Y$ to
$\psi_1(y)\, \in\, (E_G/P)_{\phi\circ\psi_1 (y)}\,=\,
((\psi^*E_G)/P)_y$.

Let
$$
E^Y_P\, \subset\, \psi^*E_G
$$
be the principal $P$--bundle over $Y$ defined by
the section $\sigma$ constructed above. The line bundle
$\psi^*_1 L_\chi$ over $Y$ is identified with $E^Y_P(\chi)$,
the line bundle associated to the principal $P$--bundle
$E^Y_P$ for the character $\chi$. Therefore, the given condition
that the degree of $E^Y_P(\chi)$ is nonnegative implies that
$\text{degree}(\psi^*_1 L_\chi)\, \geq \, 0$.
Consequently, $L_\chi$ is numerically
effective. Finally, Theorem \ref{thm.h.d.} says that $E_G$ is
semistable with $c_2(\text{ad}(E_G))\, =\, 0$. This completes the
proof of the proposition.
\end{proof}

Let $E$ be a holomorphic vector bundle of rank $r\,\geq\, 2$
over $M$. Fix an integer $k\, \in\, [1\, ,r-1]$.
Proposition \ref{pr.-1} says that $E$ is semistable
with $c_2(\text{End}(E))\, =\, 0$ if and only
if for every holomorphic map $\psi\, :\, Y\, \longrightarrow\, M$,
where $Y$ is any compact connected Riemann surface,
and for every subbundle $F\, \subset\, \psi^*E$, of rank $k$ (the
fixed integer), the degree of $F^*\otimes (\psi^*E/F)$ is nonnegative.

Let $G$ be a complex linear algebraic group such that
$G\, =\, G_1\times G_2\times \cdots \times G_{\ell}$, where
each $G_i$ is simple. A principal $G$--bundle $E_G$ over $M$
is semistable if and only if each $G_i$--bundle, $i\,\in\,
[1\, , \ell]$, obtained by extending the structure group of $E_G$
using the projection $G\, \longrightarrow\, G_i$, is semistable.
Any proper parabolic subgroup of $G$ is of the form $Q_1\times
Q_2\times \cdots \times Q_{\ell}$, where $Q_i$ is either
a parabolic subgroup of $G_i$ or $Q_i \, =\, G_i$, with
at least one $Q_j$ a parabolic subgroup.

Theorem \ref{thm.h.d.} and Proposition \ref{pr.-1} together have
the following corollary:

\begin{corollary}\label{cor.-22}
Fix a parabolic subgroup $P \,=\, Q_1\times Q_2
\times \cdots \times Q_{\ell}\, \subset\, G\,=\,
G_1\times G_2\times \cdots \times G_{\ell}$, where
$Q_i$, $i\, \in\, [1\, ,\ell]$,
is a proper parabolic subgroup of the simple linear algebraic
group $G_i$ defined over
$\mathbb C$. Fix a character $\chi$ of $P$ whose restriction
to each $Q_i$ is antidominant and nontrivial. Let $E_G$
be a principal $G$--bundle over a connected projective
manifold $M$. The following are equivalent:
\begin{enumerate}
\item{} The $G$--bundle $E_G$ is semistable and the second Chern
class
$$
c_2({\rm ad}(E_G)) \, \in\, H^4(M,\, {\mathbb Q})
$$
vanishes.

\item{} The associated line bundle $L_\chi\, :=\,
(E_G\times {\mathbb C})/P$ over $E_G/P$ for the character
$\chi$ is numerically effective.

\item{} For every pair of the form $(Y\, ,\psi)$,
where $Y$ is a compact connected Riemann surface and
$$
\psi\,:\, Y\,\longrightarrow\, M
$$
a holomorphic map, and every reduction
$E_P\, \subset\, \psi^*E_G$ of structure group to $P$ of the
$G$--bundle $\psi^*E_G$ over $Y$, the associated line bundle
$E_P(\chi)\, =\, (E_P\times {\mathbb C})/P$ over $Y$ is of
nonnegative degree.
\end{enumerate}
\end{corollary}

\begin{proof}
We have
$$
\text{ad}(E_G) \, =\, \bigoplus_{i=1}^{\ell} \text{ad}(E_{G_i})\, ,
$$
where $E_{G_i}$ is the principal $G_i$--bundle
obtained by extending the structure group of $E_G$
using the projection $G\, \longrightarrow\, G_i$. As
$c_1(\text{ad}(E_{G_i}))\, =\, 0$ for each
$i\, \in\, [1\, ,\ell]$, we have
\begin{equation}\label{eq.-su.}
c_2(\text{ad}(E_G))\,=\, \sum_{i=1}^{\ell} c_2(\text{ad}(E_{G_i}))
\,\in\, H^4(M,\, {\mathbb Q})\, .
\end{equation}
Also, as we noted earlier, $E_G$ is semistable if and only if
each $E_{G_i}$ is semistable. Therefore, Theorem \ref{thm.h.d.}
(respectively, Proposition \ref{pr.-1}) says that statement
(2) (respectively, statement (3)) implies statement (1).
Following the proof of Proposition \ref{pr.-1} we conclude 
that (2) and (3) are equivalent.

Assume that statement (1) is valid.

Since $\text{ad}(E_{G_i})$ is semistable, and
$c_1(\text{ad}(E_{G_i}))\, \in\,H^2(M,\, {\mathbb Q})$
vanishes, the Bogomolov inequality says
\begin{equation}\label{eq.-bo.}
0\, \leq\, (c_2(\text{ad}(E_{G_i}))\cup c_1(\xi)^{d-2})\cap [M]
\,\in\, {\mathbb Z}\, ,
\end{equation}
where $\xi$ is the fixed very ample line bundle
over $M$ and $d\,=\,
\dim_{\mathbb C} M$ \cite{Bo}. From Eq. \eqref{eq.-su.} we have
\begin{equation}\label{e.1}
(c_2(\text{ad}(E_{G}))\cup c_1(\xi)^{d-2})\cap [M]\,=\,
\sum_{i=1}^{\ell} (c_2(\text{ad}(E_{G_i}))\cup c_1(\xi)^{d-2})
\cap [M]\, .
\end{equation}
On the other hand, from the given condition that
$c_2(\text{ad}(E_{G}))\, =\, 0$ we know that the left--hand
side of Eq. \eqref{e.1} vanishes. Therefore, from Eq.
\eqref{eq.-bo.} it follows immediately that
\begin{equation}\label{eq.-su.2}
(c_2(\text{ad}(E_{G_i}))\cup c_1(\xi)^{d-2})\cap [M]\, =\, 0
\end{equation}
for each $i\, \in\, [1\, ,\ell]$.

As $\text{ad}(E_{G_i})$ is semistable with
$c_1(\text{ad}(E_{G_i}))\, =\, 0$ and satisfies Eq.
\eqref{eq.-su.2}, we conclude that
the vector bundle $\text{ad}(E_{G_i})$ admits a filtration
of subbundles such that each successive quotient has the property
that all the Chern classes of positive degree vanish
\cite[page 39, Theorem 2]{Si}. This immediately implies that
$c_2(\text{ad}(E_{G_i}))\, =\, 0$ for each $i\, \in\, [1\, ,\ell]$.
Consequently, using Theorem \ref{thm.h.d.} and Proposition
\ref{pr.-1} we conclude
that statement (1) implies both statements (2) and (3) (recall
that $E_G$ is semistable if and only if each $E_{G_i}$ is
semistable). This completes the proof of the corollary.
\end{proof}

Let $G$ be a connected reductive linear algebraic group
defined over $\mathbb C$.
As before, the center of $G$ will be denoted by $Z(G)$.
Therefore, $G/Z(G)$ is a product of simple groups. A parabolic
subgroup of $G$ is the inverse image, for the projection of $G$
to $G/Z(G)$, of a parabolic subgroup of $G/Z(G)$. We recall that
a parabolic subgroup $P$ of $G$ is said to be \textit{without
any simple factor} if the projection of $P/Z(G)$ to each simple
factor of $G/Z(G)$ is a proper parabolic subgroup.

A $G$--bundle $E_G$ over $M$ is semistable if and only if the
$G/Z(G)$--bundle $E_G/Z(G)$ is semistable \cite[page 146,
Proposition 7.1]{Ra}.

Therefore, using Corollary \ref{cor.-22} we have:

\begin{theorem}\label{thm.2}
Let $G$ be a connected reductive linear algebraic group
over $\mathbb C$. Fix a proper parabolic subgroup
$P\, \subset\, G$ without any simple factor, and also fix a
character $\chi$ of $P$ such that
\begin{enumerate}
\item[(i)] the character $\chi$ is trivial on the center
$Z(G)\, \subset\, G$, and
\item[(ii)] the restriction of $\chi$ to the parabolic subgroup,
defined by $P$, of each simple factor of $G/Z(G)$ is
nontrivial and antidominant.
\end{enumerate}
Let $E_G$ be a principal $G$--bundle over a connected projective
manifold $M$. Then the following four statements are equivalent:
\begin{enumerate}
\item{} The $G$--bundle $E_G$ is semistable and the
second Chern class
$$
c_2({\rm ad}(E_G)) \, \in\, H^4(M,\, {\mathbb Q})
$$
vanishes.

\item{} The associated line bundle $L_\chi\, :=\,
(E_G\times {\mathbb C}_{\chi})/P$ over $E_G/P$ for the
character $\chi$ is numerically effective.

\item{} For every pair of the form $(Y\, ,\psi)$,
where $Y$ is a compact connected Riemann surface and
$$
\psi\,:\, Y\,\longrightarrow\, M
$$
a holomorphic map, and every reduction $E_P\, \subset\,
\psi^*E_G$ of structure group to $P$ of the principal
$G$--bundle $\psi^*E_G$ over $Y$, the associated line bundle
$E_P(\chi)\, =\, (E_P\times {\mathbb C}_{\chi})/P$ over
$Y$ is of nonnegative degree.

\item{} For any pair $(Y\, ,\psi)$ as in (3), the $G$--bundle
$\psi^*E_G$ over $Y$ is semistable.
\end{enumerate}
\end{theorem}

Note that statement (3) in Theorem \ref{thm.2} implies
statement (4). Statement (4) follows from statement (3)
by using the criterion
for semistability of a $G$--bundle over a Riemann surface
given by the equivalence of statements (1) and (3).

\begin{remark}
{\rm Theorem \ref{thm.2} shows that for a principal
$G$--bundle $E_G$ over $M$, where $G$ is connected reductive,
with $c_2({\rm ad}(E_G))\, =\, 0$, the condition that
$E_G$ is semistable does not depend on the choice of the
polarization on $M$ needed for defining degree.}
\end{remark}

\section{Criterion for parabolic bundles}\label{parabolic}

Let $H$ be a connected linear
algebraic group over $\mathbb C$. Fix $n$ distinct points
$$
D \, :=\, \{x_1,\cdots , x_n\} \, \subset\, X
$$
of a compact connected Riemann surface $X$. We will refer
to the points $\{x_i\}_{i=1}^n$ as parabolic points.

Let $E'_H$ be a connected smooth quasiprojective
variety over $\mathbb C$ and
$$
f \, :\, E'_H\times H\, \longrightarrow\, E'_H
$$
an algebraic action of $H$ on $E'_H$.

A \textit{parabolic $H$--bundle} over $X$ with parabolic structure
over $D$ is a pair
$(E'_H\, ,f)$ as above together with a dominant morphism
$$
\psi\, :\, E'_H\, \longrightarrow\, X
$$
satisfying the following conditions:
\begin{enumerate}
\item{} $\psi\circ f \, =\, \psi\circ p_1$, where $p_1$ is the
projection of $E'_H\times H$ to $E'_H$, or equivalently, the
map $\psi$ is equivariant for the action of $H$;

\item{} for each point $x\, \in\, X$, the action of $H$ on the
fiber $\psi^{-1}(x)_{\text{red}}$ is transitive;

\item{} the restriction of $\psi$ to $\psi^{-1}(X\setminus D)$
makes $\psi^{-1}(X\setminus D)$ a principal $H$--bundle over
$X\setminus D$, that is, the map $\psi$ is smooth over
$\psi^{-1}(X\setminus D)$ and the map to the fiber product
$$
\psi^{-1}(X\setminus D)\times H\, \longrightarrow\,
\psi^{-1}(X\setminus D)\times_{X\setminus D} 
\psi^{-1}(X\setminus D)
$$
defined by $(z\, ,g) \longmapsto\, (z\, , f(z,g))$ is an
isomorphism;

\item{} for each point $z\, \in\,
\psi^{-1}(D)_{\text{red}}$, the isotropy
at $z$ for the action of $H$ is a finite subgroup of $H$.
\end{enumerate}
See \cite{BBN} for the details. See \cite{MS} for parabolic
vector bundles.

For notational convenience, a parabolic $H$--bundle defined
as above will be denoted by $E_*$.

Given a parabolic $\text{GL}(n,{\mathbb C})$--bundle $E_*$ over
$X$, using the standard action of $\text{GL}(n,{\mathbb C})$ on
${\mathbb C}^n$ the principal $\text{GL}(n,{\mathbb C})$--bundle
over $X\setminus D$ defined by $E_*$ gives a vector bundle
over $X\setminus D$. This vector bundle has a natural extension
to $X$,
which is constructed using $E_*$, that carries the parabolic
structure of the parabolic vector bundle corresponding to $E_*$.

Let $Y$ be a compact connected Riemann surface and
$\Gamma\, \subset\, \text{Aut}(Y)$ a finite subgroup
of the group of all holomorphic automorphisms of $Y$. Let
$E_H$ be a principal $H$--bundle over $Y$. A
\textit{$\Gamma$--linearization} of $E_H$ is a lift of the
action of $\Gamma$ on $Y$ to the total space of $E_H$ that
commutes with the action of $H$. So a $\Gamma$--linearization
of $E_H$ is a left action of $\Gamma$ on $E_H$ such
that for any $\gamma\, \in\, \Gamma$, the automorphism
of the variety $E_H$ defined by it is an isomorphism of the
$H$--bundle $E_H$ over the automorphism $\gamma$ of $Y$.

Let $E_*$ be a parabolic $H$--bundle over $X$ with $D$ as the
parabolic divisor. There is a (ramified) finite
Galois covering
\begin{equation}\label{gal.c.}
\phi\, :Y\, \longrightarrow\, X
\end{equation}
and a $\Gamma$--linearized $H$--bundle $E_H$ over $Y$,
where $\Gamma$ is the Galois group of the covering $\phi$,
such that $E_H$ corresponds to $E_*$ (see \cite{BBN}). The
covering $\phi$ is ramified over $D$ and for any $x\, \in\,
D$, the order of ramification is a multiple of the order of the
isotropy subgroup for any point in $\psi^{-1}(x)_{\text{red}}$.
See \cite[Ch. 1.1, pages 303--305]{KMM} for the construction of
such a covering.

A reduction of structure group to a closed subgroup
$P\, \subset\, H$ of a parabolic
$H$--bundle $E_*$ is a holomorphic section
$$
\sigma\, :\, X\, \longrightarrow\, E'_H/P
$$
of the natural projection of $E'_H/P$ to $X$, where $E'_H$ is
the underlying variety for $E_*$. Note that $q^{-1}(\sigma(X))$,
where $q\, :\, E'_H \,\longrightarrow\,
E'_H/P$ is the natural projection, is a
parabolic $P$--bundle. Conversely, if $E'_P\, \subset\, E'_H$
is a parabolic $P$--bundle with the
induced action of $P$, then it defines a section
$\sigma$ as above. This section $\sigma$ has the property that
$q^{-1}(\sigma(X))$ coincides with $E'_P$.

If $E_H$ is a $\Gamma$--linearized
principal $H$--bundle over $Y$ corresponding
to the parabolic $H$--bundle $E_*$ over $X$, then reductions
of $E_*$ to $P$ are in bijective correspondence with the
$\Gamma$--invariant reductions of $E_H$ to $P$ (see \cite{BBN}).

Let $G$ be a connected
reductive linear algebraic group over $\mathbb C$.
Let $E_*$ be a parabolic $G$--bundle with $E'_G$ as the
underlying $G$--variety. The
parabolic $G$ bundle $E_*$ is called \textit{semistable}
if for any reduction of structure group
$E'_P\, \subset\, E'_G$ to any parabolic
subgroup $P\, \subset\, G$ and for any antidominant
character $\chi$ of $P$ trivial on the
center of $G$ the parabolic degree of the
associated parabolic line bundle $E'_P(\chi)$ is nonnegative.

If $E_G$ is a $\Gamma$--linearized
principal $G$--bundle over $Y$ corresponding
to the parabolic $G$--bundle $E_*$ over $X$, then
$E_*$ is semistable if and only if $E_G$ is so
(see \cite{BBN}).

For a parabolic $G$--bundle $E_*$ with $E'_G$ as the
underlying variety, let $N\, =\,
N(E_*)$ be the least common multiple of the order
of the isotropy groups for the action of $G$ on
$E'_G$.

Let $E'_P\, \subset\, E'_G$ be a reduction of structure
group to $P$. If $\chi$ is a character of $P$, then
$\chi^{N}$ is trivial on all
the isotropy subgroups for the action of $P$ on $E'_P$
(as $N$ is a multiple of the order of each isotropy
subgroup).
Therefore, the associated line bundle $E'_P(\chi^{N})$
has trivial parabolic structure.

Let $E_*$ be a parabolic $G$--bundle, where $G$ is a simple
linear algebraic group,
over $X$ with $E'_G$ as the underlying variety.
If $\chi$ is a character of a parabolic subgroup
$P\, \subset\, G$, then we noted that
$\chi^{N}$ is trivial on all
the isotropy subgroups for the action of $P$ on $E'_G$,
where, as before, $N\, =\, N(E_*)$. Therefore, the quotient
$$
E'_G(\chi^N) \, :=\, (E'_G\times {\mathbb C})/P
$$
is a line bundle over $E'_G/P$, with $P$ acting on
${\mathbb C}$ through $\chi^N$.

Theorem \ref{thm.2} gives
the following criteria for semistability of $E_*$.

\begin{proposition}\label{prop.-1}
Fix a pair $(P\, ,\chi)$, where $P\, \subset\, G$ is a proper
parabolic subgroup of the reductive group
and $\chi\, :\, P\, \longrightarrow\,
{\mathbb C}^*$ a nontrivial antidominant character. The parabolic
$G$--bundle $E_*$ is semistable if and only if the associated
line bundle $E'_G(\chi^N)$ over $E'_G/P$
is numerically effective, where $N$ is the least common
multiple of the order of the isotropy subgroups for the action
of $G$ on $E'_G$.

The parabolic $G$--bundle
$E_*$ is semistable if and only if for every
(possibly ramified) covering $\psi\, :\, Y\, \longrightarrow\, X$,
and for every reduction of structure group $E'_P\, \subset\,
\psi^*E'_G$ to $P$ of the pulled back parabolic
$G$--bundle $\psi^*E'_G\, =\, Y\times_X E'_G$, the associated
line bundle $E_P(\chi^N)$ over $Y$ is of nonnegative degree.
\end{proposition}

\medskip
\noindent
\textbf{Acknowledgements.} The first author thanks the International
Centre for Theoretical Physics, Trieste, for hospitality, as well as
the International Mathematical Union for a travel
support. The second author's research is
partly supported by the National Research Project ``Geometria
delle variet\`a algebriche''. The second author is a member
of the research group {\sc vbac} (Vector Bundles on Algebraic
Curves), which is partially supported by {\sc eager} ({\sc ec fp5}
Contract no. {\sc hprn-ct-2000-00099}) and
by {\sc edge} ({\sc ec fp5} Contract no. {\sc hprn-ct-2000-00101}).


\end{document}